# Strapdown Attitude Computation: Functional Iterative Integration versus Taylor Series Expansion

Yuanxin Wu and Yury A. Litmanovich

*Abstract*— This paper compares two basic approaches to solving ordinary differential equations, which form the basis for attitude computation in strapdown inertial navigation systems, namely, the Taylor series expansion approach that was used in its low-order form for deriving all mainstream algorithms and the functional iterative integration approach developed recently. They are respectively applied to solve the kinematic equations of major attitude parameters, including the quaternion, the Rodrigues vector and the rotation vector. Specifically, the mainstream algorithms, which have relied on the simplified rotation vector without exception, are considerably extended by the Taylor series expansion approach using the exact rotation vector and recursive calculation of high-order derivatives. The functional iterative integration approach is respectively implemented on both the normal polynomial and the Chebyshev polynomial. Numerical results under the classical coning motion are reported to assess all derived attitude algorithms. It is revealed that in the relative frequency range when the coning to sampling frequency ratio is below 0.05-0.1 (depending on the chosen polynomial truncation order), all algorithms have the same order of accuracy if the same number of samples are used to fit the angular velocity over the iteration interval; in the range of higher relative frequency, the group of Quat/Rod/RotFIter algorithms (by the functional iterative integration approach combined with the Chebyshev polynomial) perform the best in both accuracy and robustness, thanks to the excellent numerical stability and powerful functional representation capability of the Chebyshev polynomial.

## I. INTRODUCTION

Attitude information is vitally important for moving objects in many areas including unmanned vehicle navigation and control, virtual/augmented reality, satellite communication, robotics, and computer vision [1]. Integrating gyroscope-measured angular velocity information is an essential and self-contained way to acquire the attitude, rotation or orientation [2-5]. A number of orientation parameters have been used for attitude computation, including but not limited to the Euler angle, the rotation vector, the direction cosine matrix and the quaternion. Attitude computation is in essence to numerically solve the ordinary differential

The first author of the paper was supported in part by National Key R&D Program of China (2018YFB1305103) and National Natural Science Foundation of China (61673263). A short version will be presented at IEEE/ION Position Location and Navigation Symposium (PLANS) in 2020.

Authors' address: Yuanxin Wu, Shanghai Key Laboratory of Navigation and Location-based Services, School of Electronic Information and Electrical Engineering, Shanghai Jiao Tong University, Shanghai, China, 200240, E-mail: (yuanx_wu@hotmail.com); Yury A. Litmanovich, Central Scientific and Research Institute "Elektropribor", Saint Petersburg, Russia, 197046, E-mail: (yulitman@mail.ru).

equations of these attitude parameters. In the early years of strapdown inertial navigation systems, Savage [6] tried the Picard-type iterative technique to integrate the direction cosine rate equation, and the NASA technical report [7] made the angular velocity polynomial approximation from a sequence of gyroscope outputs and then integrated the direction cosine matrix rate by the Runge-Kutta method. Shortly after in the 1970s, the modern-day strapdown attitude algorithm structure was established on the Taylor series expansion approach by Jordan and Bortz [8, 9], which, without exception, has relied on the approximate rotation vector for incremental attitude update [10-17]. In parallel, a number of related fields employ the quaternion to deal with attitude computation, e.g., robotics [18, 19], space applications [20] and computational mathematics [21-23], where the structure-preserving attribute of geometric integration is mostly concerned.

It has long been believed that the modern-day attitude algorithm is already good enough for applications [2, 13, 15]. However, the present-day dynamic applications and the future precision gyroscopes [24, 25] demand more accurate attitude algorithms, in view of the pratically limited gyroscope output rate and the fundamental approximation of modern-day algorithms. In principle, the recent advances in [26-32] have shown that higher attitude accuracy can be achieved by dealing with better or exact rate forms of attitude parameters and using high-order numerical integration methods [33-35]. For instance, the impact of commonly-neglected third term in the rotation rate vector was evaluated in [36] and it was partially incorporated in the algorithms to gain significant accuracy in dynamic applications [27, 37]. The Rodrigues vector and the quaternion were used for attitude computation by way of functional iterative integration and Chebyshev polynomial approximation [29-31]. The functional iterative integration approach combined with Chebyshev polynomial approximation was developed independently in the navigation community [29-31, 38-41], and lately found to closely resemble the so-called Picard-Chebyshev method that was dated back to as early as 1960s [42]. In the 1980s-1990s and even quite recently, it was employed and advanced by researchers in the field of astrodynamics for orbital determination [43-46]. The Picard iterative technique actually originated from the Picard–Lindelöf theorem that has been commonly used to prove the existence and uniqueness of numerical solutions of ordinary differential equations [47], but it has been hardly exploited in engineering applications because the repeated computation of integrals are often conceived to be inconvenient and tedious [48]. The Taylor series method experiences a similar story [48] in that it is conceptually easy to work with but the high-order derivatives are taken as being tedious and complicated to calculate. To avoid the need for high-order derivatives, the Runge-Kutta methods were thus devised while attempting to retain the accuracy of the Taylor series approximation [33]. The seminal paper by Miller [10] used a low-order Taylor series with low-order angular velocity approximation to solve the approximate rotation vector rate equation. Very recently, the Taylor series approach is employed to directly solve the direction cosine matrix rate equation [28], kind of a retrospective work into the early attempt [7].

Another trend of intensive investigation in the strapdown attitude research was the ad-hoc algorithm optimization by means of special tuning of the algorithm coefficients to reduce the coning drift under the assumed motions such as the classical/generalized coning motion, the regular precession and the stochastic angular motion [10-12, 17, 36, 49, 50]. The cost one should pay is that the optimized algorithms rank below the corresponding original algorithms in accuracy in practically irregular angular motions, say maneuvers [16]. Needless to say it is desirable for algorithms to exhibit the same order of accuracy regardless of input attitude motions.

It should be noted that the angular velocity polynomial approximation from a sequence of gyroscope outputs is an integral part of both approaches [7, 10, 16, 33, 34, 36, 37], explicitly or implicitly. Hereby in this paper we consider zero and single integral of angular velocity as gyroscope measurements because they are common in real systems, although multiple integrals could also be accounted for as done in [36, 37, 51]. The purpose of this paper is to investigate the problem of strapdown attitude computation from the general perspective of solving the differential equations involved, so as to make a comprehensive approach comparison from the accuracy standpoint. It was motivated by a long fruitful discussion of the two authors about the actual superiority of the functional iterative approach over the traditional Taylor expansion approach.

The technical contribution of the paper is multiple-folded. Specifically, the Taylor-based method is substantially extended using the exact rotation vector and recursive calculation of high-order derivatives; the functional iteration method is further implemented by using the normal polynomial; the approach superiority in the previous works is refined as an immediate result of the extended Taylor-based method. The remainder of the paper is organized as follows: Section II briefly reviews two basic approaches to solve the ordinary differential equations, namely, the Taylor series expansion and the functional iterative integration. Section III discusses angular velocity approximation by the normal polynomial and the Chebyshev polynomial, and then Sections IV-VI make use of the two basic approaches to solve the kinematic equations of major attitude parameters including the quaternion, the Rodrigues vector and the rotation vector. Specifically, Sections V-VI both rely on the functional iterative integration, yet with the normal polynomial and the Chebyshev polynomial, respectively. Section VII comprehensively assesses all derived algorithms under the classical coning motion. A brief summary is given in the last section of the paper.

## II. General Approaches to Solve Ordinary Differential Equation

Without the loss of generality, consider an ordinary differential equation over a time interval $\begin{bmatrix} 0 & t \end{bmatrix}$

$$\dot{\mathbf{y}} = \mathbf{f}\left(\mathbf{y}, t\right) \tag{1}$$

where $\mathbf{f}\left(\mathbf{y}, t\right)$ is an infinitely smooth function and the initial value of $\mathbf{y}$ is given by $\mathbf{y}\left(0\right)$.

The solution to (1) could be obtained by the Taylor series expanded at $t = 0$

$$\mathbf{y}\left(t\right) = \mathbf{y}\left(0\right) + \dot{\mathbf{y}}\left(0\right)t + \ddot{\mathbf{y}}\left(0\right)\frac{t^2}{2!} + \cdots = \sum_{j=0}^{\infty} \mathbf{y}^{(j)}\left(0\right)\frac{t^j}{j!} \tag{2}$$

in which $\mathbf{y}^{(j)}(0) = \mathbf{y}^{(j)}(t)\big|_{t=0}$, that is, the value of the $i$-th derivative at $t = 0$. The series is infinite and keeping terms up to the $m$-th order leads to a Taylor series approximation as follows [33]

$$\mathbf{y}(t) = \sum_{j=0}^{m} \mathbf{y}^{(j)}(0)\frac{t^j}{j!} + \mathbf{y}^{(m+1)}(\varsigma)\frac{t^{m+1}}{(m+1)!}, \quad \varsigma \in \begin{bmatrix} 0 & t \end{bmatrix} \tag{3}$$

The second term on the right side characterizes the error of the $m$-th order Taylor series approximation. It is commonly conceived that the calculation of high-order derivatives involved is time-consuming and tedious, though conceptually straightforward [33]. In fact, the high-order derivatives share the common property that only their values at $t = 0$ are required and thus we need not know the analytic forms of the derivatives. The values of the high-order derivatives could be recursively computed by making use of the calculus rule of elementary functions [34, 48]. For example, assume $\mathbf{f}(\mathbf{y}, t) = \mathbf{y}(t)\mathbf{u}(t)$, then $\mathbf{y}^{(j)}(0) = \left(\mathbf{y}(t)\mathbf{u}(t)\right)^{(j-1)}\big|_{t=0} = \left(\dot{\mathbf{y}}(t)\mathbf{u}(t) + \mathbf{y}(t)\dot{\mathbf{u}}(t)\right)^{(j-2)}\big|_{t=0} = \ldots = \sum_{i=0}^{j-1}\binom{j-1}{i}\mathbf{y}^{(i)}(0)\mathbf{u}^{(j-1-i)}(0)$, where $\binom{m}{k}$ denotes the number of combinations taking $k$ elements from $m$ elements. It means that high-order derivatives at some instant can be represented by low-order derivatives at the same instant, which provides an economical way to compute the Taylor series approximation [34, 48]. Furthermore, the solution to (1) could alternatively be obtained by the Picard iteration or a kind of functional iterative integration as [48]

$$\mathbf{y}_j(t) = \mathbf{y}(0) + \int_0^t \mathbf{f}\left(\mathbf{y}_{j-1}(t), t\right)dt, \quad j = 1, 2, \ldots \tag{4}$$

where the initial function over the integration interval could be set to $\mathbf{y}_0(t) = \mathbf{y}(0)$. It can be proved that the difference between the $(m\text{-}1)$-th and the $m$-th iterations [48]

$$\left\|\mathbf{y}_m(t) - \mathbf{y}_{m-1}(t)\right\| \le W L^{m-1} \frac{t^m}{m!} \tag{5}$$

if the function $\mathbf{f}$ is bounded by $W$, namely $W = \max_{\tau \in [0\,t]}\left\|\mathbf{f}(\mathbf{y}, \tau)\right\|$, and satisfies the Lipschitz continuity condition $\left\|\mathbf{f}(\mathbf{y}, t) - \mathbf{f}(\mathbf{z}, t)\right\| \le L\left\|\mathbf{y} - \mathbf{z}\right\|$. The right side of (5) is a term of the Taylor series for $e^{Lt}$ up to a scale. By the Weierstrass M-Test [52], the above sequence $\left\{\mathbf{y}_j\right\}_{j=0}^{\infty}$ converges uniformly to the true solution. The practical drawback of (4) is the repeated computation of integrals. It was surmounted in practice, for the first time to our best knowledge, by the Chebyshev polynomial approximation of the function $\mathbf{f}$ in [42]. Of course, the normal polynomial could alternatively be used.

It should be highlighted hereby that the Taylor series approximation in (3) is in itself of the normal polynomial, while the Picard iteration in (4) naturally accommodates any kind of polynomial.

## III. Angular Velocity Polynomial Fitted from Gyroscope Measurements

The above two basic approaches in last section requires $\mathbf{f}(\mathbf{y}, t)$ to be analytically known, but for the attitude computation under investigation, only the equally-spaced discrete gyroscope measurements are available. A common practice is to approximate the angular velocity by a polynomial fitted from the discrete gyroscope measurements [7, 10, 16, 33, 34, 36, 37]. Note that the content of this section has been presented in previous works [29-31] and is repeated here for easy reference.

Assume discrete measurements (or called samples) of angular velocity $\boldsymbol{\omega}_k$ or angular increment $\Delta\boldsymbol{\theta}_k$ are available by a triad of gyroscopes at time instants $t_k = kT$ ( $k = 1, 2, \dots N$ ), where $T$ denotes the sampling interval.

### A. Normal Polynomial

The angular velocity can be approximated by a normal polynomial as

$$\boldsymbol{\omega}(t) = \sum_{i=0}^{n} \mathbf{d}_i t^i, \ \ n \le N-1 \tag{6}$$

where the coefficient $\mathbf{d}_i$ is determined using the discrete angular velocity or angular increment measurements. The derivatives of the fitted angular velocity can be readily obtained as $\boldsymbol{\omega}^{(j)}(0) = j! \mathbf{d}_j$ for $j \le n$ and $\boldsymbol{\omega}^{(j)}(0) = 0$ otherwise.

For the case of angular velocity measurement, the coefficients $\mathbf{d}_i$ satisfy the equation

$$\begin{bmatrix} 1 & t_1 & \dots & t_1^n \\ 1 & t_2 & \dots & t_2^n \\ \vdots & \vdots & \vdots & \vdots \\ 1 & t_N & \dots & t_N^n \end{bmatrix} \begin{bmatrix} \mathbf{d}_0^T \\ \mathbf{d}_1^T \\ \vdots \\ \mathbf{d}_n^T \end{bmatrix} = \begin{bmatrix} \boldsymbol{\omega}_1^T \\ \boldsymbol{\omega}_2^T \\ \vdots \\ \boldsymbol{\omega}_N^T \end{bmatrix} \tag{7}$$

And for the case of angular increment measurement instead, the coefficients satisfy

$$\begin{bmatrix} t_1 & \dfrac{t_1^2}{2} & \dots & \dfrac{t_1^{n+1}}{n+1} \\ t_2 - t_1 & \dfrac{t_2^2 - t_1^2}{2} & \dots & \dfrac{t_2^{n+1} - t_1^{n+1}}{n+1} \\ \vdots & \vdots & \vdots & \vdots \\ t_N - t_{N-1} & \dfrac{t_N^2 - t_{N-1}^2}{2} & \dots & \dfrac{t_N^{n+1} - t_{N-1}^{n+1}}{n+1} \end{bmatrix} \begin{bmatrix} \mathbf{d}_0^T \\ \mathbf{d}_1^T \\ \vdots \\ \mathbf{d}_n^T \end{bmatrix} = \begin{bmatrix} \Delta\boldsymbol{\theta}_1^T \\ \Delta\boldsymbol{\theta}_2^T \\ \vdots \\ \Delta\boldsymbol{\theta}_N^T \end{bmatrix} \tag{8}$$

### B. Chebyshev Polynomial

The Chebyshev polynomial is a sequence of orthogonal polynomial bases and has better numerical stability than the normal polynomial [35]. The Chebyshev polynomial of the first kind is defined over the interval $\begin{bmatrix} -1 & 1 \end{bmatrix}$ by the recurrence relation as

$$F_0(x) = 1, \ F_1(x) = x, \ F_{i+1}(x) = 2xF_i(x) - F_{i-1}(x) \tag{9}$$

where $F_i(x)$ is the $i$th-degree Chebyshev polynomial of the first kind. For any $j, k \geq 0$, the Chebyshev polynomial of first kind satisfies the equality [35]

$$F_j(\tau) F_k(\tau) = \frac{1}{2}\left(F_{j+k}(\tau) + F_{|j-k|}(\tau)\right). \tag{10}$$

In order to apply the Chebyshev polynomial, the actual time interval $\begin{bmatrix} 0 & t_N \end{bmatrix}$ is mapped onto $\begin{bmatrix} -1 & 1 \end{bmatrix}$ by letting $t = (1+\tau) t_N/2$. Then, the angular velocity over the mapped interval is fitted by the Chebyshev polynomial, given by

$$\boldsymbol{\omega}(\tau) = \sum_{i=0}^{n} \mathbf{c}_i F_i(\tau), \ n \leq N-1 \tag{11}$$

The coefficient $\mathbf{c}_i$ is determined for the case of angular velocity measurement by solving the equation as follows:

$$\begin{bmatrix} 1 & F_1(\tau_1) & \cdots & F_n(\tau_1) \\ 1 & F_1(\tau_2) & \cdots & F_n(\tau_2) \\ \vdots & \vdots & \vdots & \vdots \\ 1 & F_1(\tau_N) & \cdots & F_n(\tau_N) \end{bmatrix} \begin{bmatrix} \mathbf{c}_0^T \\ \mathbf{c}_1^T \\ \vdots \\ \mathbf{c}_n^T \end{bmatrix} = \begin{bmatrix} \boldsymbol{\omega}_1^T \\ \boldsymbol{\omega}_2^T \\ \vdots \\ \boldsymbol{\omega}_N^T \end{bmatrix} \tag{12}$$

According to the integral property of the Chebyshev polynomial [35], we have

$$G_{i,[\tau_{k-1} \tau_k]} \triangleq \int_{\tau_{k-1}}^{\tau_k} F_i(\tau) d\tau = \begin{cases} \left(\dfrac{iF_{i+1}(\tau_k)}{i^2-1} - \dfrac{\tau_k F_i(\tau_k)}{i-1}\right) - \left(\dfrac{iF_{i+1}(\tau_{k-1})}{i^2-1} - \dfrac{\tau_{k-1} F_i(\tau_{k-1})}{i-1}\right), & i \neq 1 \\ \dfrac{\tau_k^2 - \tau_{k-1}^2}{2}, & i = 1 \end{cases} \tag{13}$$

With the aid of $(11)$ and $(13)$, the angular increment is related to the fitted angular velocity by

$$\Delta\boldsymbol{\theta}_{t_k} = \int_{t_{k-1}}^{t_k} \boldsymbol{\omega} \, dt \overset{t \to \tau}{=} \frac{t_N}{2} \int_{\tau_{k-1}}^{\tau_k} \boldsymbol{\omega} \, d\tau = \frac{t_N}{2} \sum_{i=0}^{n} \mathbf{c}_i \int_{\tau_{k-1}}^{\tau_k} F_i(\tau) d\tau = \frac{t_N}{2} \sum_{i=0}^{n} \mathbf{c}_i G_{i,[\tau_{k-1} \tau_k]} \tag{14}$$

Then, the coefficient $\mathbf{c}_i$ in $(11)$ is determined for the case of angular increment measurement by solving the following equation:

$$\begin{bmatrix} G_{0,[\tau_0 \tau_1]} & G_{1,[\tau_0 \tau_1]} & \cdots & G_{n,[\tau_0 \tau_1]} \\ G_{0,[\tau_1 \tau_2]} & G_{1,[\tau_1 \tau_2]} & \cdots & G_{n,[\tau_1 \tau_2]} \\ \vdots & \vdots & \vdots & \vdots \\ G_{0,[\tau_{N-1} \tau_N]} & G_{1,[\tau_{N-1} \tau_N]} & \cdots & G_{n,[\tau_{N-1} \tau_N]} \end{bmatrix} \begin{bmatrix} \mathbf{c}_0^T \\ \mathbf{c}_1^T \\ \vdots \\ \mathbf{c}_n^T \end{bmatrix} = \frac{2}{t_N} \begin{bmatrix} \Delta\boldsymbol{\theta}_1^T \\ \Delta\boldsymbol{\theta}_2^T \\ \vdots \\ \Delta\boldsymbol{\theta}_N^T \end{bmatrix} \tag{15}$$

The linear equations $(7)$, $(8)$, $(12)$ and $(15)$ could be well solved by the common least-square method.

## IV. Attitude Algorithms by Taylor Series Expansion

This section will apply the Taylor series expansion to solve the kinematic equations of quaternion, Rodrigues vector and rotation vector for attitude computation, by employing recursive calculation of high-order derivatives.

*A. QuatTaylor by Quaternion*

The attitude quaternion kinematic equation is related to the angular velocity as [4]

$$\dot{\mathbf{q}} = \frac{1}{2}\mathbf{q} \circ \boldsymbol{\omega} \tag{16}$$

Attitude quaternion $\mathbf{q}$ is represented as a four-dimensional column vector of unit magnitude, i.e., $\mathbf{q} = \begin{bmatrix} s & \boldsymbol{\eta}^T \end{bmatrix}^T$, where $s$ is the scalar part and $\boldsymbol{\eta}$ is the vector part. If the scalar and vector parts are regarded as a scalar quaternion and a vector quaternion, respectively, then quaternion can be alternatively written as $\mathbf{q} = s + \boldsymbol{\eta}$. The product of two quaternions is given by

$\mathbf{q}_1 \circ \mathbf{q}_2 = \begin{bmatrix} s_1 & -\boldsymbol{\eta}_1^T \\ \boldsymbol{\eta}_1 & s_1 \mathbf{I}_3 + \boldsymbol{\eta}_1 \times \end{bmatrix} \begin{bmatrix} s_2 \\ \boldsymbol{\eta}_2 \end{bmatrix} = \begin{bmatrix} s_1 s_2 - \boldsymbol{\eta}_1^T \boldsymbol{\eta}_2 \\ s_1 \boldsymbol{\eta}_2 + s_2 \boldsymbol{\eta}_1 + \boldsymbol{\eta}_1 \times \boldsymbol{\eta}_2 \end{bmatrix}$. Denote the Euler rotation axis and the Euler rotation angle by $\mathbf{e}$ and

$\alpha$ respectively, the quaternion can be alternatively expressed as $\mathbf{q} = \cos\dfrac{\alpha}{2} + \mathbf{e}\sin\dfrac{\alpha}{2}$. $\boldsymbol{\omega}$ is the angular velocity vector quaternion with zero scalar part, formed by the three-dimensional angular velocity vector.

The $j$-th order derivative of the quaternion can be recursively computed as

$$\mathbf{q}^{(j)}(0) = \frac{1}{2}\left(\mathbf{q}\circ\boldsymbol{\omega}\right)^{(j-1)}\Big|_{t=0} = \frac{1}{2}\sum_{i=0}^{j-1}\binom{i-1}{j}\mathbf{q}^{(j-1-i)}(0)\circ\boldsymbol{\omega}^{(i)}(0) \tag{17}$$

Explicitly, $\mathbf{q}^{(0)}(0) = \mathbf{q}(0)$, $\mathbf{q}^{(1)}(0) = \dfrac{1}{2}\mathbf{q}(0)\circ\mathbf{d}_0$ and $\mathbf{q}^{(2)}(0) = \dfrac{1}{2}\left(\mathbf{q}^{(1)}(0)\circ\mathbf{d}_0 + \mathbf{q}^{(0)}(0)\circ\mathbf{d}_1\right)$, etc.

Then, the Taylor series approximation (3) can be explicitly written as

$$\mathbf{q}(t) = \frac{1}{2}\sum_{j=0}^{m}\left(\sum_{i=0}^{j-1}\binom{i-1}{j}\mathbf{q}^{(j-1-i)}(0)\circ\boldsymbol{\omega}^{(i)}(0)\right)\frac{t^j}{j!} \tag{18}$$

*B. RodTaylor by Rodrigues Vector*

The Rodrigues vector rate equation is related to the angular velocity as [1, 16]

$$\dot{\mathbf{g}} = \boldsymbol{\omega} + \frac{1}{2}\mathbf{g}\times\boldsymbol{\omega} + \frac{1}{4}\mathbf{g}\mathbf{g}^T\boldsymbol{\omega} \tag{19}$$

The Rodrigues vector, $\mathbf{g} = 2\tan(\alpha/2)\mathbf{e}$, is transformed to the attitude quaternion by $\mathbf{q} = \dfrac{2+\mathbf{g}}{\sqrt{4+|\mathbf{g}|^2}}$. Note that it is now used for incremental attitude update, i.e., $\mathbf{g}(0) = 0$. The $j$-th order derivative of the Rodrigues vector can be recursively computed as

$$\mathbf{g}^{(j)}(0) = \left( \boldsymbol{\omega} + \frac{1}{2}\mathbf{g} \times \boldsymbol{\omega} + \frac{1}{4}\mathbf{g}\mathbf{g}^T\boldsymbol{\omega} \right)^{(j-1)} \Bigg|_{t=0}$$

$$= \boldsymbol{\omega}^{(j-1)}(0) + \frac{1}{2}\sum_{i=0}^{j-1}\binom{j-1}{i}\mathbf{g}^{(j-1-i)}(0) \times \boldsymbol{\omega}^{(i)}(0) + \frac{1}{4}\sum_{i=0}^{j-1}\binom{j-1}{i}\mathbf{g}^{(j-1-i)}(0)\sum_{k=0}^{i}\binom{i}{k}\left(\mathbf{g}^{(i-k)}(0)\right)^T\boldsymbol{\omega}^{(k)}(0)$$

$$= \boldsymbol{\omega}^{(j-1)}(0) + \frac{1}{2}\sum_{i=0}^{j-1}\binom{j-1}{i}\mathbf{g}^{(j-1-i)}(0) \times \boldsymbol{\omega}^{(i)}(0) + \frac{1}{4}\sum_{i=0}^{j-1}\sum_{k=0}^{i}\binom{j-1}{i}\binom{i}{k}\mathbf{g}^{(j-1-i)}(0)\left(\mathbf{g}^{(i-k)}(0)\right)^T\boldsymbol{\omega}^{(k)}(0) \tag{20}$$

Substituting into (3) yields the Taylor series approximation of the Rodrigues vector

$$\mathbf{g}(t) = \sum_{j=0}^{m}\mathbf{g}^{(j)}(0)\frac{t^j}{j!} \tag{21}$$

*C. RotTaylor by Rotation Vector*

The rotation vector's rate equation is related to the angular velocity by [9]

$$\dot{\boldsymbol{\sigma}} = \boldsymbol{\omega} + \frac{1}{2}\boldsymbol{\sigma} \times \boldsymbol{\omega} + \frac{1}{|\boldsymbol{\sigma}|^2}\left(1 - \frac{|\boldsymbol{\sigma}|\sin|\boldsymbol{\sigma}|}{2(1-\cos|\boldsymbol{\sigma}|)}\right)\boldsymbol{\sigma} \times (\boldsymbol{\sigma} \times \boldsymbol{\omega}) \triangleq \boldsymbol{\omega} + \frac{1}{2}\boldsymbol{\sigma} \times \boldsymbol{\omega} + A(|\boldsymbol{\sigma}|)\boldsymbol{\sigma} \times (\boldsymbol{\sigma} \times \boldsymbol{\omega}) \tag{22}$$

The rotation vector, $\boldsymbol{\sigma} = \alpha\mathbf{e}$, is transformed to the attitude quaternion by $\mathbf{q} = \cos\frac{|\boldsymbol{\sigma}|}{2} + \frac{\boldsymbol{\sigma}}{|\boldsymbol{\sigma}|}\sin\frac{|\boldsymbol{\sigma}|}{2}$ for nonzero $\boldsymbol{\sigma}$, typically used for incremental attitude update as well, i.e., $\boldsymbol{\sigma}(0) = 0$. The scalar coefficient $A(|\boldsymbol{\sigma}|)$ is a trigonometric function of the rotation vector's magnitude and is singular at zero $\boldsymbol{\sigma}$. The $A$'s odd-order derivatives are all zeros and the leading even-order derivatives at $t = 0$ (equivalently at zero $\boldsymbol{\sigma}$) can be readily obtained from its Taylor series expansion as

$$A(|\boldsymbol{\sigma}|) = \frac{1}{12} + \frac{|\boldsymbol{\sigma}|^2}{720} + \frac{|\boldsymbol{\sigma}|^4}{30240} + \frac{|\boldsymbol{\sigma}|^6}{1209600} + \frac{|\boldsymbol{\sigma}|^8}{47900160} + \frac{691|\boldsymbol{\sigma}|^{10}}{1307674368000} + \cdots \tag{23}$$

$$A(0) = 1/12, A^{(2)}(0) = 1/360, A^{(4)}(0) = 1/1260$$
$$A^{(6)}(0) = 1/1680, A^{(8)}(0) = 1/1188, A^{(10)}(0) = 691/360360, \ldots \tag{24}$$

Using (22), the *j*-th order derivative of the rotation vector can be recursively computed as

$$\boldsymbol{\sigma}^{(j)}(0) = \left( \boldsymbol{\omega} + \frac{1}{2}\boldsymbol{\sigma} \times \boldsymbol{\omega} + A\boldsymbol{\sigma} \times (\boldsymbol{\sigma} \times \boldsymbol{\omega}) \right)^{(j-1)} \Bigg|_{t=0}$$

$$= \boldsymbol{\omega}^{(j-1)}(0) + \frac{1}{2}\sum_{i=0}^{j-1}\binom{j-1}{i}\boldsymbol{\sigma}^{(j-1-i)}(0) \times \boldsymbol{\omega}^{(i)}(0) + \sum_{i=0}^{j-1}\binom{j-1}{i}A^{(j-1-i)}(0)\sum_{k=0}^{i}\binom{i}{k}\boldsymbol{\sigma}^{(i-k)}(0) \times \sum_{s=0}^{k}\binom{k}{s}C_k^s\left(\boldsymbol{\sigma}^{(k-s)}(0) \times \boldsymbol{\omega}^{(s)}(0)\right)$$

$$= \boldsymbol{\omega}^{(j-1)}(0) + \frac{1}{2}\sum_{i=0}^{j-1}\binom{j-1}{i}\boldsymbol{\sigma}^{(j-1-i)}(0) \times \boldsymbol{\omega}^{(i)}(0) + \sum_{i=0}^{j-1}\sum_{k=0}^{i}\sum_{s=0}^{k}\binom{j-1}{i}\binom{i}{k}\binom{k}{s}A^{(j-1-i)}(0)\boldsymbol{\sigma}^{(i-k)}(0) \times \left(\boldsymbol{\sigma}^{(k-s)}(0) \times \boldsymbol{\omega}^{(s)}(0)\right) \tag{25}$$

Substituting into (3) yields the Taylor series approximation of the rotation vector

$$\boldsymbol{\sigma}(t) = \sum_{j=0}^{m}\boldsymbol{\sigma}^{(j)}(0)\frac{t^j}{j!} \tag{26}$$

The algorithm presented in [36] is a special case of RotTaylor for $N = 3$ that exactly considers up to the fifth-order derivatives of the rotation vector. Of special interest are two approximate rotation vectors, overwhelmingly used by the mainstream attitude algorithms [10-15] and related to the angular velocity by

$$\dot{\boldsymbol{\sigma}} \approx \boldsymbol{\omega} + \frac{1}{2}\boldsymbol{\sigma} \times \boldsymbol{\omega} \quad \text{or} \quad \dot{\boldsymbol{\sigma}} \approx \boldsymbol{\omega} + \frac{1}{2}\left(\int_0^t \boldsymbol{\omega}\, dt\right) \times \boldsymbol{\omega} \tag{27}$$

For the above approximate rotation vectors, the $j$-th order derivative of the first one is simply the sum of the first two terms in (25), i.e.,

$$\boldsymbol{\sigma}^{(j)}(0) \approx \left(\boldsymbol{\omega} + \frac{1}{2}\boldsymbol{\sigma} \times \boldsymbol{\omega}\right)^{(j-1)}\bigg|_{t=0} = \boldsymbol{\omega}^{(j-1)}(0) + \frac{1}{2}\sum_{i=0}^{j-1}\binom{j-1}{i}\boldsymbol{\sigma}^{(j-1-i)}(0) \times \boldsymbol{\omega}^{(i)}(0) \tag{28}$$

And the $j$-th order derivative of the second one is

$$\boldsymbol{\sigma}^{(j)}(0) \approx \left(\boldsymbol{\omega} + \frac{1}{2}\left(\int_0^t \boldsymbol{\omega}\, dt\right) \times \boldsymbol{\omega}\right)^{(j-1)}\bigg|_{t=0} = \boldsymbol{\omega}^{(j-1)}(0) + \frac{1}{2}\sum_{i=0}^{j-1}\binom{j-1}{i}\boldsymbol{\omega}^{(j-2-i)}(0) \times \boldsymbol{\omega}^{(i)}(0) \tag{29}$$

Thus, the Taylor series approximation of the rotation vector (26), with derivatives given by (28) or (29), is particularly named as RotTaylor-T2 or RotTaylor-T2s (hereby 's' means further simplification). These approximations are of special interest in that RotTaylor-T2/T2s encompass almost all of the mainstream attitude algorithms in the literature that are founded on the approximate rotation vectors (27). In particular, the relationship of RotTaylor-T2s to the mainstream 2/3-sample algorithms [10] is shown in Appendix.

It should be noted that the Taylor series approximations in (18), (21) and (26) are time polynomials that actually reconstruct the whole attitude history over the time interval $\begin{bmatrix} t_0 & t_N \end{bmatrix}$ in which the $N$ gyroscope samples are measured, sharing the same advantage of the attitude algorithms by functional iterative integration [16, 33, 34].

## V. ATTITUDE ALGORITHMS BY FUNCTIONAL ITERATIVE INTEGRATION (NORMAL POLYNOMIAL)

This section will use the functional iterative integration approach to solve the differential equations of attitude parameters, using the normal polynomial approximation. Notably, the functional iterative integration, combined with the Chebyshev polynomial approximation, has been successfully applied for attitude computation in [29-31]. By analogy, the development procedure is straightforward.

*A. QuatFIter-np[1] by Quaternion*

With the angular velocity polynomial in (6), the functional iterative integration is applied to the attitude quaternion rate equation (16), yielding

$$\mathbf{q}_j(t) = \mathbf{q}(0) + \frac{1}{2}\int_0^t \mathbf{q}_{j-1} \circ \boldsymbol{\omega} \, dt = \mathbf{q}(0) + \frac{1}{2}\int_0^t \mathbf{q}_{j-1}(t) \circ \sum_{i=0}^n \mathbf{d}_i t^i \, dt, \quad j = 1, 2, \ldots \tag{30}$$

Suppose the attitude quaternion at the $(j\text{-}1)$-th iteration is represented by a normal polynomial of order $m_{j-1}$, i.e.,

$\mathbf{q}_{j-1}(t) = \sum_{k=0}^{m_{j-1}} \mathbf{b}_{j-1,k} t^k$. Substituting into (30) gives

$$\begin{aligned}
\mathbf{q}_j(t) &= \mathbf{q}(0) + \frac{1}{2}\int_0^t \sum_{k=0}^{m_{j-1}} \mathbf{b}_{j-1,k} t^k \circ \sum_{i=0}^n \mathbf{d}_i t^i \, dt \\
&= \mathbf{q}(t_0) + \frac{1}{2}\sum_{k=0}^{m_{j-1}} \sum_{i=0}^n \frac{\mathbf{b}_{j-1,k} \circ \mathbf{d}_i}{k+i+1} t^{k+i+1}, \quad j = 1, 2, \ldots
\end{aligned} \tag{31}$$

It can be seen that the polynomial order of attitude quaternion grows quickly by $m_j = m_{j-1} + n + 1$. In the explicit form,

$\mathbf{q}_0(t) = \mathbf{q}(0)$ and $\mathbf{q}_1(t) = \mathbf{q}(0) + \frac{1}{2}\mathbf{q}(0) \circ \sum_{i=0}^n \frac{\mathbf{d}_i}{i+1} t^{i+1}$, etc.

*B. RodFIter-np by Rodrigues Vector*

Apply the functional iterative integration to the Rodrigues vector's rate equation (19),

$$\begin{aligned}
\mathbf{g}_j(t) &= \int_0^t \left( \boldsymbol{\omega} + \frac{1}{2}\mathbf{g}_{j-1} \times \boldsymbol{\omega} + \frac{1}{4}\mathbf{g}_{j-1}\mathbf{g}_{j-1}^T \boldsymbol{\omega} \right) dt \\
&= \int_0^t \boldsymbol{\omega} \, dt + \frac{1}{2}\int_0^t \mathbf{g}_{j-1} \times \boldsymbol{\omega} \, dt + \frac{1}{4}\int_0^t \mathbf{g}_{j-1}\mathbf{g}_{j-1}^T \boldsymbol{\omega} \, dt, \quad j = 1, 2, \ldots
\end{aligned} \tag{32}$$

Suppose the Rodrigues vector at the $(j\text{-}1)$-th iteration is represented by a normal polynomial of order $m_{j-1}$, i.e.,

$\mathbf{g}_{j-1}(t) = \sum_{k=0}^{m_{j-1}} \mathbf{b}_{j-1,k} t^k$. Using the angular velocity polynomial approximation (6), the three integrals on the right side of (32) can

be computed as

$$\int_0^t \boldsymbol{\omega} \, dt = \int_0^t \sum_{i=0}^n \mathbf{d}_i t^i \, dt = \sum_{i=0}^n \frac{\mathbf{d}_i}{i+1} t^{i+1} \tag{33}$$

$$\int_0^t \mathbf{g}_{j-1} \times \boldsymbol{\omega} \, dt = \int_0^t \sum_{k=0}^{m_{j-1}} \mathbf{b}_{j-1,k} t^k \times \sum_{i=0}^n \mathbf{d}_i t^i \, dt = \sum_{k=0}^{m_{j-1}} \sum_{i=0}^n \frac{\mathbf{b}_{j-1,k} \times \mathbf{d}_i}{k+i+1} t^{k+i+1} \tag{34}$$

---

[1] Abbreviation 'FIter' stands for Functional Iterative integration; 'np' stands for normal polynomial.

$$\int_0^t \mathbf{g}_{j-1} \mathbf{g}_{j-1}^T \boldsymbol{\omega} \, dt = \int_0^t \sum_{s=0}^{m_{j-1}} \mathbf{b}_{j-1,s} t^s \sum_{k=0}^{m_{j-1}} \mathbf{b}_{j-1,k}^T t^k \sum_{i=0}^n \mathbf{d}_i t^i \, dt = \sum_{s=0}^{m_{j-1}} \sum_{k=0}^{m_{j-1}} \sum_{i=0}^n \frac{\mathbf{b}_{j-1,s} \mathbf{b}_{j-1,k}^T \mathbf{d}_i}{s+k+i+1} t^{s+k+i+1}$$

(35)

Substituting (33)-(35) into (32), the Rodrigues vector at the $j$-th iteration is obtained as

$$\mathbf{g}_j(t) = \sum_{i=0}^n \frac{\mathbf{d}_i}{i+1} t^{i+1} + \frac{1}{2} \sum_{k=0}^{m_{j-1}} \sum_{i=0}^n \frac{\mathbf{b}_{j-1,k} \times \mathbf{d}_i}{k+i+1} t^{k+i+1} + \frac{1}{4} \sum_{s=0}^{m_{j-1}} \sum_{k=0}^{m_{j-1}} \sum_{i=0}^n \frac{\mathbf{b}_{j-1,s} \mathbf{b}_{j-1,k}^T \mathbf{d}_i}{s+k+i+1} t^{s+k+i+1}$$

(36)

Obviously, the polynomial order grows by $m_j = 2m_{j-1} + n + 1$.

### C. RotFIter-np by Rotation Vector

The rotation vector's rate equation (22) involves trigonometric functions and it is cumbersome to do the integrals. Here the coefficient of the third term is approximated by $A \approx 1/12$.

Suppose the Rodrigues vector at the $(j\text{-}1)$-th iteration is represented by a normal polynomial of order $m_{j-1}$, i.e.,

$\boldsymbol{\sigma}_{j-1}(t) = \sum_{k=0}^{m_{j-1}} \mathbf{b}_{j-1,k} t^k$. By analogy with the development in (33)-(36), the rotation vector at the $j$-th iteration can be written as

$$\boldsymbol{\sigma}_j(t) \approx \sum_{i=0}^n \frac{\mathbf{d}_i}{i+1} t^{i+1} + \frac{1}{2} \sum_{k=0}^{m_{j-1}} \sum_{i=0}^n \frac{\mathbf{b}_{j-1,k} \times \mathbf{d}_i}{k+i+1} t^{k+i+1} + \frac{1}{12} \sum_{s=0}^{m_{j-1}} \sum_{k=0}^{m_{j-1}} \sum_{i=0}^n \frac{\mathbf{b}_{j-1,s} \times (\mathbf{b}_{j-1,k} \times \mathbf{d}_i)}{s+k+i+1} t^{s+k+i+1}$$

(37)

The polynomial order grows by $m_j = 2m_{j-1} + n + 1$ as well. It is particularly named as RotFIter-np-T3, because the third term is only approximately accounted for. For the case of only considering the first two terms as in (27), the corresponding algorithm is readily obtained by abandoning the last additive term in (37), named as RotFIter-np-T2.

Note that Eqs. (31), (36) and (37) could be iterated by updating the normal polynomial coefficients only and truncating the normal polynomials at each iteration to avoid fast order growing, as done in QuatFIter [33]. The polynomial truncation order, denoted by $m_r$ hereafter, also acts as the highest order of derivative to those algorithms by the Taylor series expansion. The iteration times could be controlled by some pre-defined maximum or stopping criterion given in the sequel.

### VI. ATTITUDE ALGORITHMS BY FUNCTIONAL ITERATIVE INTEGRATION (CHEBYSHEV POLYNOMIAL)

This section is mainly a brief summary of the proposed algorithms in [29-31] for easy reference and comparison in this paper. Readers are referred to those works for details. It should be highlighted that the technique has been successfully applied to the whole inertial navigation algorithm including velocity/position computation [40, 41].

### A. QuatFIter by Quaternion

With the angular velocity polynomial given by (11), the functional iterative integration approach is applied to solve the quaternion rate equation [31]

$$\mathbf{q}_j = \mathbf{q}(0) + \frac{1}{2}\int_0^t \mathbf{q}_{j-1} \circ \boldsymbol{\omega}\, dt = \mathbf{q}(0) + \frac{t_N}{4}\int_{-1}^\tau \mathbf{q}_{j-1} \circ \boldsymbol{\omega}\, d\tau. \tag{38}$$

Assume $\boldsymbol{\omega}$ and the quaternion estimate at the $(j\text{-}1)$-th iteration is given by a weighted sum of Chebyshev polynomials, say

$\mathbf{q}_j \triangleq \sum_{k=0}^{m_{j-1}} \mathbf{b}_{j-1,k} F_k(\tau)$ where $m_{j-1}$ is the maximum degree and $\mathbf{b}_{j-1,k}$ is the coefficient of the $k^{\text{th}}$-degree Chebyshev polynomial

at the $(j\text{-}1)$-th iteration. Substituting it, together with (10) and (11), into (38)

$$\mathbf{q}_j(\tau) = \mathbf{q}(0) + \frac{t_N}{8}\sum_{k=0}^{m_{j-1}}\sum_{i=0}^{n} \mathbf{b}_{j-1,k} \circ \mathbf{c}_i \left( G_{k+i,[-1\ \tau]} + G_{|k-i|,[-1\ \tau]} \right) \tag{39}$$

where $G_{i,[-1\ \tau]}$ is the integrated $i^{\text{th}}$-degree Chebyshev polynomial over the interval $\begin{bmatrix} -1 & \tau \end{bmatrix}$, as defined in (13).

### B. RodFIter by Rodrigues Vector

When the functional iterative integration approach is applied to solve the Rodrigues vector rate equation, we get [29, 30]

$$\mathbf{g}_j = \int_0^t \left( \mathbf{I}_3 + \frac{1}{2}\mathbf{g}_{j-1}\times + \frac{1}{4}\mathbf{g}_{j-1}\mathbf{g}_{j-1}^T \right) \boldsymbol{\omega}\, dt = \frac{t_N}{2}\int_{-1}^\tau \left( \mathbf{I}_3 + \frac{1}{2}\mathbf{g}_{j-1}\times + \frac{1}{4}\mathbf{g}_{j-1}\mathbf{g}_{j-1}^T \right) \boldsymbol{\omega}\, d\tau \tag{40}$$

Assume the Rodrigues vector at the $(j\text{-}1)$-th iteration is given by a weighted sum of Chebyshev polynomials, say

$\mathbf{g}_{j-1} \triangleq \sum_{k=0}^{m_{j-1}} \mathbf{b}_{j-1,k} F_k(\tau)$. Substituting into (40) leads to

$$\mathbf{g}_j(\tau) = \frac{t_N}{2}\left( \begin{array}{l} \sum_{i=0}^n \mathbf{c}_i G_{i,[-1\ \tau]} + \frac{1}{4}\sum_{k=0}^{m_{j-1}}\sum_{i=0}^n \mathbf{b}_{j-1,k}\times\mathbf{c}_i \left( G_{k+i,[-1\ \tau]} + G_{|k-i|,[-1\ \tau]} \right) \\ + \frac{1}{16}\sum_{s=0}^{m_{j-1}}\sum_{k=0}^{m_{j-1}}\sum_{i=0}^n \mathbf{b}_{j-1,s}\mathbf{b}_{j-1,k}^T\mathbf{c}_i \left( G_{s+k+i,[-1\ \tau]} + G_{|s+k-i|,[-1\ \tau]} + G_{|s-k|+i,[-1\ \tau]} + G_{||s-k|-i|,[-1\ \tau]} \right) \end{array} \right) \tag{41}$$

### C. RotFIter by Rotation Vector

The functional iterative integration approach is applied to solve the full rate equation of the rotation vector

$$\boldsymbol{\sigma}_j = \int_0^t \left( \mathbf{I}_3 + \frac{1}{2}\boldsymbol{\sigma}_{j-1}\times + A(|\boldsymbol{\sigma}|)(\boldsymbol{\sigma}_{j-1}\times)^2 \right)\boldsymbol{\omega}\, dt = \frac{t_N}{2}\int_{-1}^\tau \left( \mathbf{I}_3 + \frac{1}{2}\boldsymbol{\sigma}_{j-1}\times + A(|\boldsymbol{\sigma}|)(\boldsymbol{\sigma}_{j-1}\times)^2 \right)\boldsymbol{\omega}\ d\tau \tag{42}$$

Suppose the Rodrigues vector at the $(j\text{-}1)$-th iteration is represented by a Chebyshev polynomial of order $m_{j-1}$, i.e.,

$\boldsymbol{\sigma}_{j-1}(t) = \sum_{k=0}^{m_{j-1}} \mathbf{b}_{j-1,k} F_k(\tau)$. Recall that the scalar coefficient $A(|\boldsymbol{\sigma}|)$ involves a trigonometric function of the rotation vector's magnitude and does not allow a tractable analytical integration.

Denote the third term by $A(|\boldsymbol{\sigma}|)\boldsymbol{\sigma}_{j-1}\times(\boldsymbol{\sigma}_{j-1}\times\boldsymbol{\omega}) \triangleq \boldsymbol{\eta}_{j-1}$ and approximate it by a Chebyshev polynomial, say $\boldsymbol{\eta}_{j-1} = \sum_{k=0}^{p_{j-1}} \boldsymbol{\gamma}_{j-1,k} F_k(\tau)$.

The coefficients can be approximately computed as [35]

$$\boldsymbol{\gamma}_{j-1,k} \approx \frac{2-\delta_{0k}}{Q} \sum_{s=0}^{Q-1} \cos\left(\frac{k(s+1/2)\pi}{Q}\right) \boldsymbol{\eta}_{j-1}\left(\cos\left(\frac{(s+1/2)\pi}{Q}\right)\right) \tag{43}$$

where $\delta_{ij}$ is the Kronecker delta function. Exact coefficients could be obtained only if the number of summation terms $Q$ approaches infinity. Substituting into (42) gives

$$\boldsymbol{\sigma}_j(\tau) = \frac{t_N}{2}\left(\sum_{i=0}^{n}\mathbf{c}_i G_{i,[-1\ \tau]} + \frac{1}{4}\sum_{k=0}^{m_{j-1}}\sum_{i=0}^{n}\mathbf{b}_{j-1,k}\times\mathbf{c}_i\left(G_{k+i,[-1\ \tau]} + G_{|k-i|,[-1\ \tau]}\right) + \sum_{k=0}^{p_{j-1}}\boldsymbol{\gamma}_{j-1,k}G_{k,[-1\ \tau]}\right) \tag{44}$$

When $A(|\boldsymbol{\sigma}|) \approx 1/12$, the above algorithm reduces to the RodFIter-T3 given in [16], i.e.,

$$\boldsymbol{\sigma}_j(\tau) = \frac{t_N}{2}\left(\begin{array}{c}\sum_{i=0}^{n}\mathbf{c}_i G_{i,[-1\ \tau]} + \frac{1}{4}\sum_{k=0}^{m_{j-1}}\sum_{i=0}^{n}\mathbf{b}_{j-1,k}\times\mathbf{c}_i\left(G_{k+i,[-1\ \tau]} + G_{|k-i|,[-1\ \tau]}\right) \\ + \frac{1}{48}\sum_{s=0}^{m_{j-1}}\sum_{k=0}^{m_{j-1}}\sum_{i=0}^{n}\mathbf{b}_{j-1,s}\times\mathbf{b}_{j-1,k}\times\mathbf{c}_i\left(\begin{array}{c}G_{s+k+i,[-1\ \tau]} + G_{|s+k-i|,[-1\ \tau]} \\ +G_{|s-k|+i,[-1\ \tau]} + G_{||s-k|-i|,[-1\ \tau]}\end{array}\right)\end{array}\right) \tag{45}$$

Totally omitting the third term of the rotation vector rate equation (42) gives us the RotFIter-T2 in [16].

Similiar to the last section, Eqs. (39), (41) and (44) could be iterated by updating the Chebyshev polynomial coefficients only and making necessary polynomial truncation at each iteration to avoid fast order growing [16, 33, 34].

## VII. NUMERICAL RESULTS AND ALGORITHM COMPARISON

All derived attitude algorithms in the paper can be directly implemented by following the given formulae, as summarized in Table I. Simulations are performed in this section under the classical coning motion scenario to evaluate these algorithms. The coning motion has explicit analytical expressions in the angular velocity and the attitude parameter, so it has been widely accepted as a standard criterion for algorithm accuracy assessment in the inertial navigation field [2, 4, 7]. It is not uncommon in practice with a large excitation of attitude drift error, e.g., in situations of angular vibration or complex rotation. The angular velocity of the classical coning motion is described by $\boldsymbol{\omega} = \Omega\left[-2\sin^2(\alpha/2) \quad -\sin(\alpha)\sin(\Omega t) \quad \sin(\alpha)\cos(\Omega t)\right]^T$, with the true rotation vector $\boldsymbol{\sigma} = \alpha\left[0 \quad \cos(\Omega t) \quad \sin(\Omega t)\right]^T$, the true Rodrigues vector $\mathbf{g} = 2\tan(\alpha/2)\left[0 \quad \cos(\Omega t) \quad \sin(\Omega t)\right]^T$ and the true quaternion $\mathbf{q} = \cos(\alpha/2) + \sin(\alpha/2)\left[0 \quad \cos(\Omega t) \quad \sin(\Omega t)\right]^T$. In the above, $\alpha$ denotes the coning angle and $\Omega = 2\pi f_c$ denotes the angular frequency of the coning motion (unit: rad/s) and $f_c$ is the coning frequency (unit: Hz). The angular increment measurement is assumed and the sampling rate is nominally set to $f_s = 1000$ Hz.

The following principal angle metric is used to quantify the attitude computation error

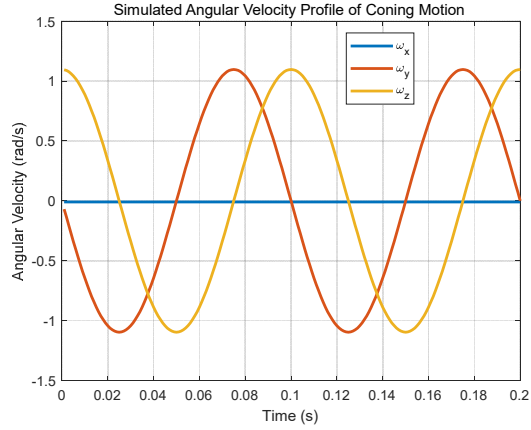

Figure 1. Angular velocity profile of classical coning motion.

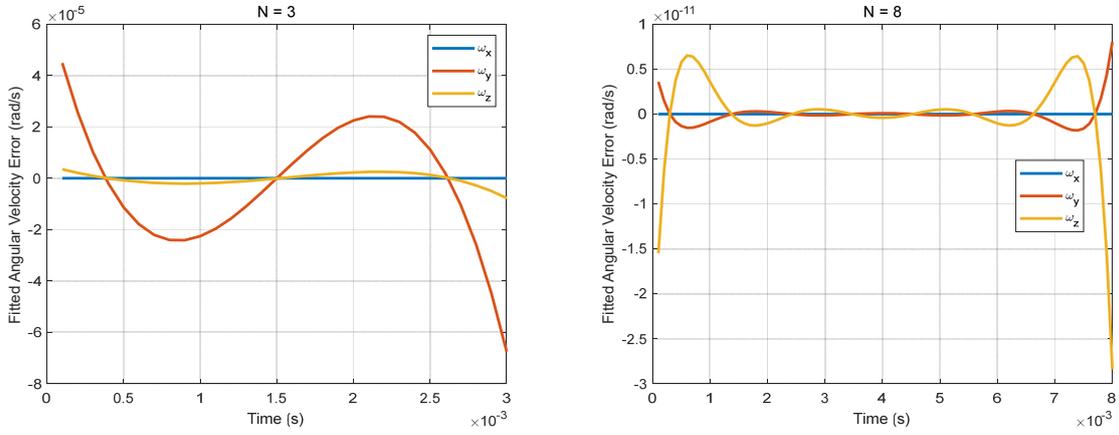

Figure 2. Errors of fitted angular velocity for $N = 3$ and 8.

$$\varepsilon_{att} = 2 \left| \left[ \mathbf{q}^* \circ \hat{\mathbf{q}} \right]_{2:4} \right| \tag{46}$$

where $\hat{\mathbf{q}}$ denotes the quaternion estimate computed by attitude algorithms, and the operator $\left[ \cdot \right]_{2:4}$ takes the vector part of the error quaternion. If the used attitude parameter is other than quaternion, then the computed result needs to be transformed to the corresponding quaternion for error quantification. The polynomial order of the fitted angular velocity (6) is uniformly set to $n = N - 1$, if not explicitly stated.

*A.  Fitted Angular Velocity Polynomial and Reconstructed Attitude*

Figure 1 plots the angular velocity profile of the classical coning motion, where the coning angle is set to $\alpha = 1$ degree with the coning frequency $f_c = 10$ Hz. Figure 2 presents the errors of the fitted angular velocity by the normal polynomial (6) or by the Chebyshev polynomial (11) during the first update interval, for the number of samples $N = 3$ and 8. It shows that using more samples leads to much more accurate fitted angular velocity, making it possible to acquire more accurate attitude. Note that an erroneous angular velocity cannot in general be compensated in the subsequent attitude computation process. The work of RodFIter [16] has thrown lights on this fact in the case of the Rodrigues vector (Theorem 2 therein). The two kinds of polynomials have identical angular velocity fitting errors but their coefficients differ much (cf. Fig. 4 in the sequel).

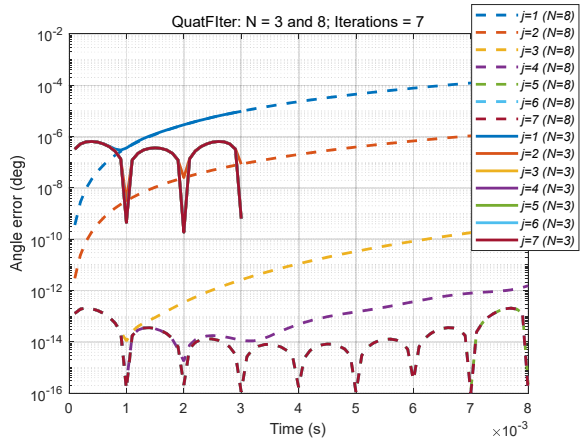

Figure 3. Attitude errors of QuatFIter across seven iterations for $N = 3$ (sold lines) and 8 (dashed lines).

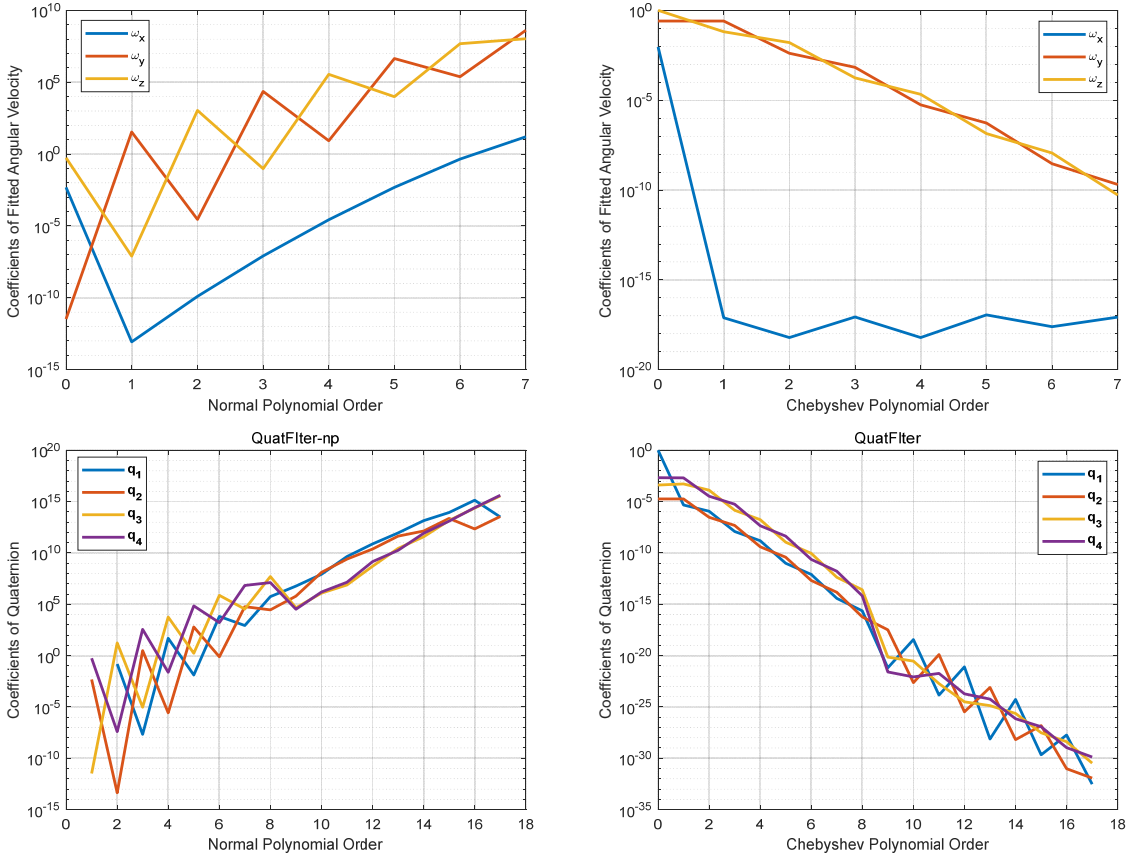

Figure 4. Polynomial coefficients of fitted angular velocity and computed quaternion by QuatFIter-np (left) and QuatFIter (right) for $N = 8$.

Taking the QuatFIter algorithm proposed in [33] as a demonstration (polynomial truncation order set to $m_T = N + 9$ ), Figure 3

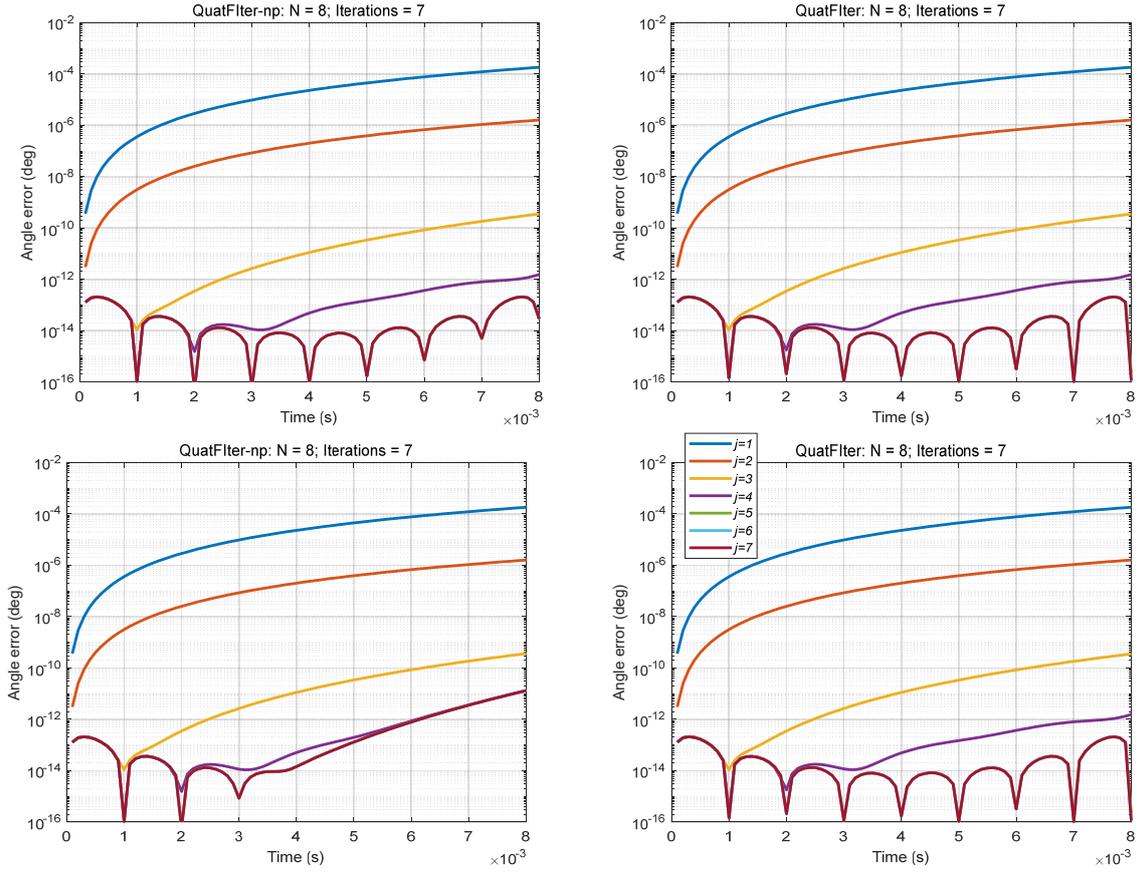

Figure 5. Truncation effect on QuatFIter-np and QuatFIter (Top: $m_T = N + 5$; bottom: $m_T = N + 2$).

plots the principal angle errors of the reconstructed attitudes over the first iteration interval, across seven iterations for the cases of $N = 3$ and 8. The angle error reduces and converges as the iteration goes on. Because the fitted angular velocity has much better accuracy, the attitude error with $N = 8$ is significantly smaller. Additionally, regarding the converged results (after two iterations with $N = 3$; after four iterations with $N = 8$), the attitude errors turn to have sharp drops at the sampling instants. This apparent 'n-shape' phenomenon is an indication of insufficient fitting of the angular velocity polynomials by the current number of gyroscope samples. Figure 4 presents the polynomial coefficients of the fitted angular velocities for the case of $N = 8$, as well as those of the computed quaternions at the 7th iteration by QuatFIter-np and QuatFIter. Along with the increasing order, the magnitude of the normal polynomial quickly increases while that of the Chebyshev polynomial swiftly decreases. The trend is observed in both the fitted angular velocities and the computed quaternions.

*B. Polynomial Truncation and Iteration Times*

Figure 5 examines the effect of two polynomial truncation orders ($m_T = N + 5, N + 2$) on QuatFIter-np and QuatFIter across seven iterations. We see that the QuatFIter-np is more vulnerable to polynomial truncation, indicating that the normal polynomial has inferior functional representation capability than the Chebyshev polynomial does. In other words, the Chebyshev polynomial requires relatively fewer terms to achieve the same accuracy, an excellent property favorable to numerical computation [35].

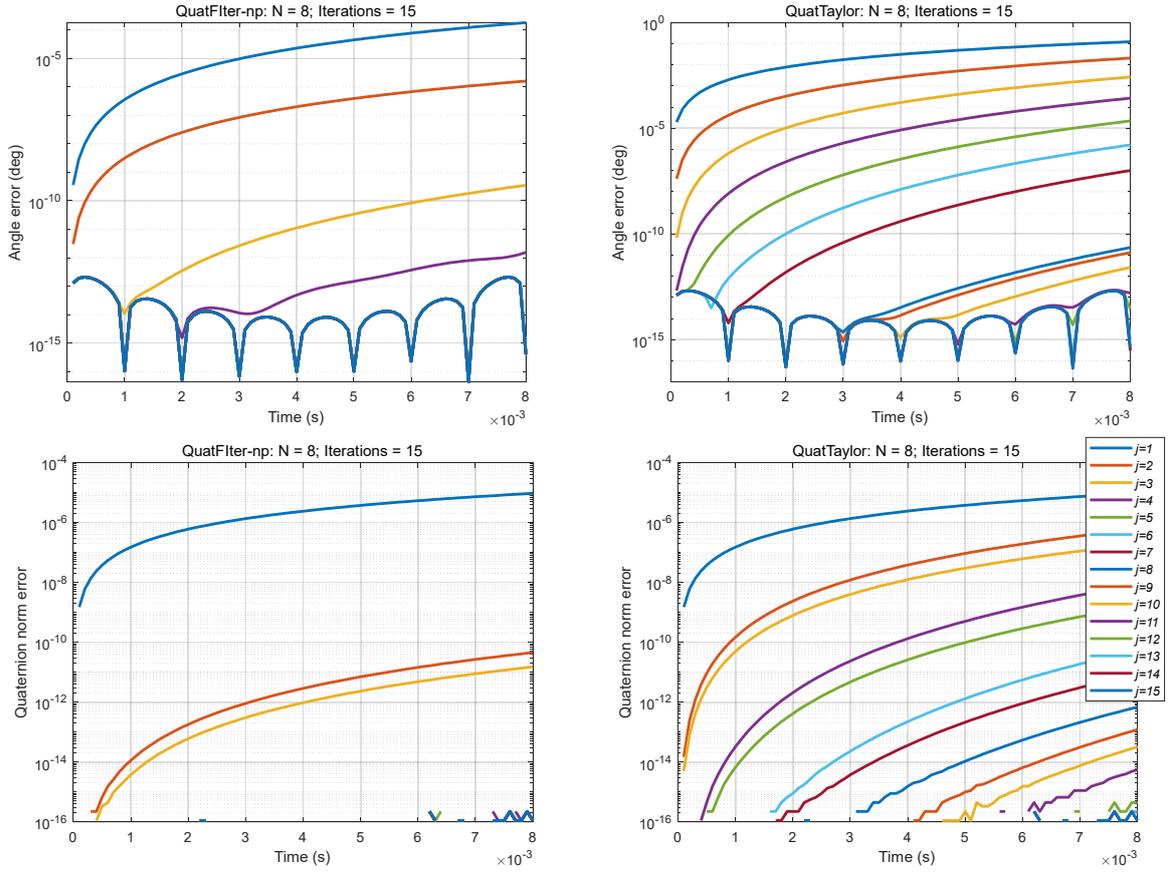

Figure 6. Attitude errors and quaternion-norm errors of QuatFIter-np and QuatTaylor. Note the order of legend colors used to denote different iterations are the same throughout the paper.

Figure 6 compares the attitude errors and quaternion-norm errors of QuatFIter-np and QuatTaylor, across fifteen iterations for the case of $N = 8$. The polynomial truncation order is still set to $m_T = N + 9$ for both algorithms. We see that it takes more iterations for QuatTaylor to reach comparable accuracy, e.g., fourteen iterations to reach the convergence in contrast to only five iterations for QuatFIter-np. This is owed to the fast increase of normal polynomial orders in QuatFIter-np, as shown in (31), while the normal polynomial order of QuatTaylor increases one by one along with each iteration. From the quaternion kinematic equation (16), it can be readily shown that the norm of quaternion is naturally preserved, i.e., $d\left(\mathbf{q}^T\mathbf{q}\right)\big/dt = 2\mathbf{q}^T\dot{\mathbf{q}} = \mathbf{q}^T\left(\mathbf{q}\circ\boldsymbol{\omega}\right) = \begin{bmatrix}1 & \mathbf{0}_{3\times1}\end{bmatrix}\boldsymbol{\omega} = 0$, if only we could compute solve the kinematic equation accurately, whether by Taylor expansion or iterative integration. We see in Fig. 6 that the quaternion norm gradually approaches to unity because the initial quaternion is of unit norm, which sets the foundation for QuatFIter [31] that uses quaternion directly for attitude computation.

It should be highlighted that the above conclusions derived from Figs. 5-6 are independent of the specific attitude parameters, which is in favor of the functional iterative integration combined with the Chebyshev polynomial for attitude computation. For a stopping criterion of iteration, the algorithms by Taylor series expansion could check if the highest-order term (HOT) is negligible relative to the require attitude accuracy, and those algorithms by functional iterative integration could use the discrepancy of polynomial coefficients (DPC) between successive iterations, namely, $\sqrt{\sum_{k=0}^{m_T} \left| \mathbf{b}_{j+1,k} - \mathbf{b}_{j,k} \right|^2}$ .

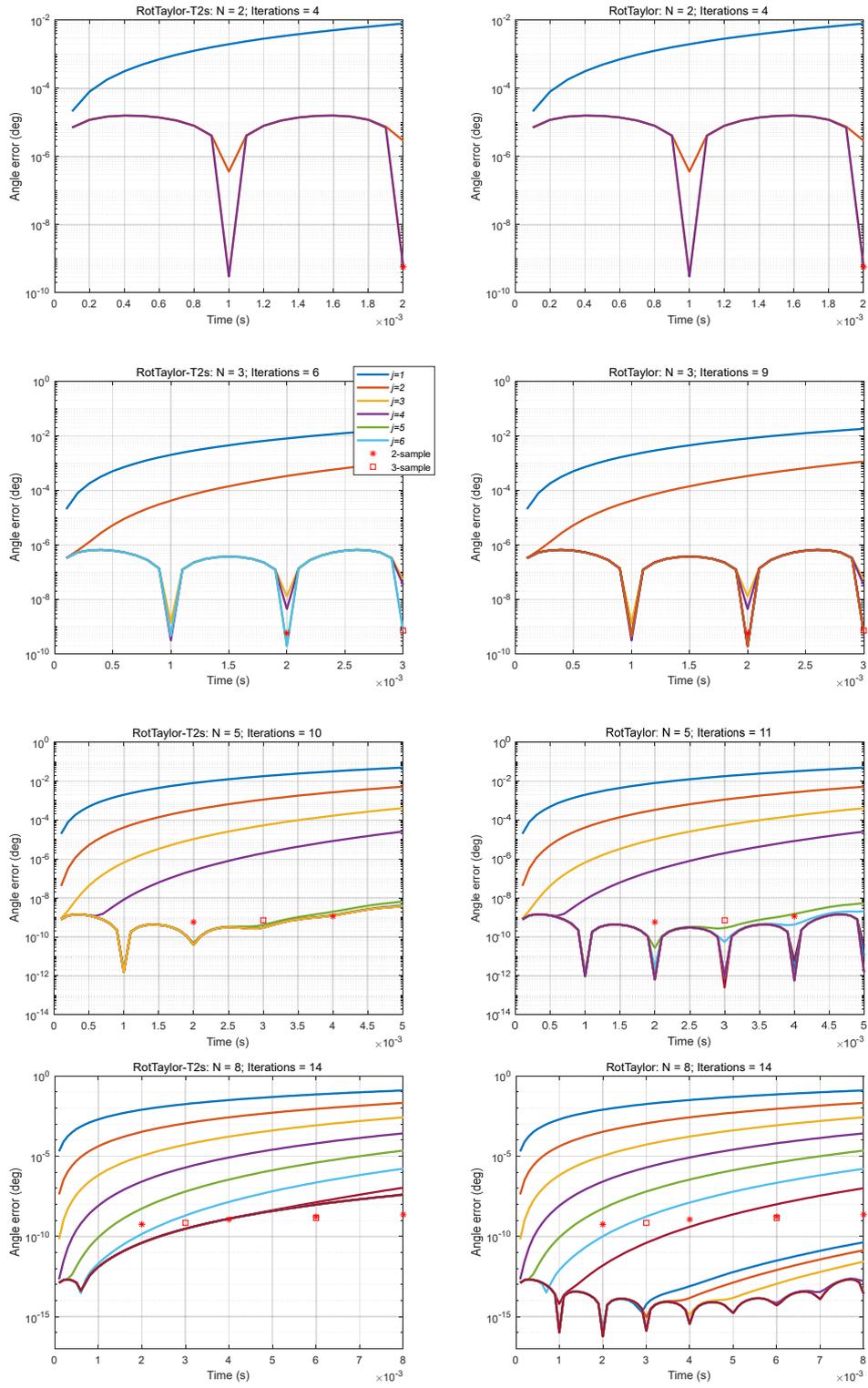

Figure 7. Attitude errors of RotTaylor-T2s (left column) and RotTaylor (right column) for $N$ = 2, 3, 5 and 8, as compared with mainstream 2-sample and 3-sample algorithms.

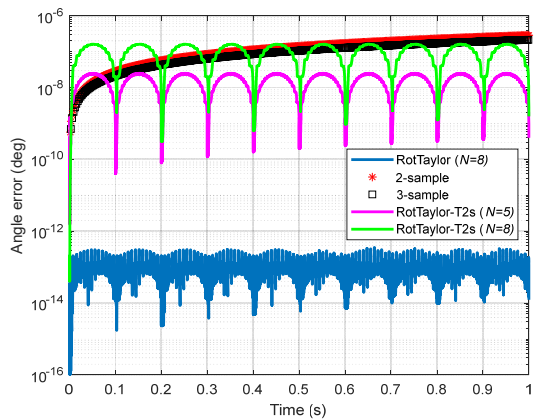

Figure 8. Attitude error comparison of RotTaylor-T2s and RotTaylor for $N$ = 5 and 8 in one second, as compared with mainstream 2/3-sample algorithms.

### C. *RotTaylor and Its Relation to Mainstream Algorithms*

Figure 7 presents the attitude errors of RotTaylor-T2s and RotTaylor for cases of $N$ = 2, 3, 5 and 8, as compared with the widely-used mainstream 2-sample and non-coning-optimized 3-sample algorithms [10, 14]. The HOT iteration-stopping criterion is employed, with truncation order $m_T = N + 9$. As the RotTaylor-T2 algorithm is founded on the approximation rotation vector, using more samples even leads to larger computation errors, as already clarified in several previous works, e.g., [33, 53]. The underlying reason is that along with more accurate angular velocity (by more samples), the RotTaylor-T2s just converges to a fake 'rotation vector' whose rate equation is exactly represented by (27) (see Theorem 2 in [16]). Only when the exact rotation vector rate equation is properly handled, e.g., by the RotTaylor just derived in this paper, could an improved accuracy be really acquired. Their attitude accuracies in one second are further demonstrated in Fig. 8.

### D. *Accuracy Comparison*

A comprehensive accuracy comparison is performed next for the case of $N$ = 8, with the coning frequency ranging 1-200 Hz. The HOT and DPC iteration-stopping criteria are respectively used in the Taylor expansion-derived algorithms and the functional iterative integration-derived algorithms. Figure 9 presents the attitude errors accumulated over one second (the analogy of the attitude error drift [36]) as a function of relative frequency $f_c / f_s$ for QuatFIter-np (QuatTaylor), RodFIter-np (RodTaylor), RotTaylor and QuatFIter (RodFIter and RotFIter), against the mainstream 2/3-sample algorithms. The results of QuatTaylor,

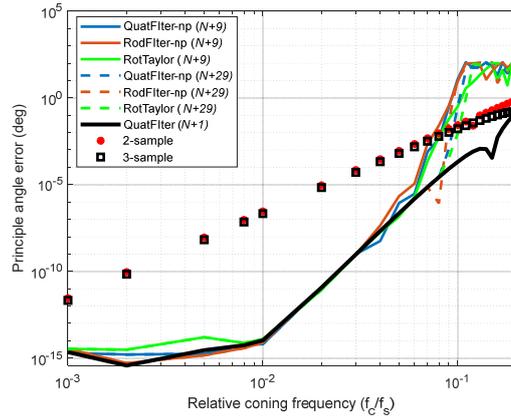

Figure 9. Attitude errors as function of relative frequency for $N = 8$, as compared with the mainstream 2/3-sample algorithms.

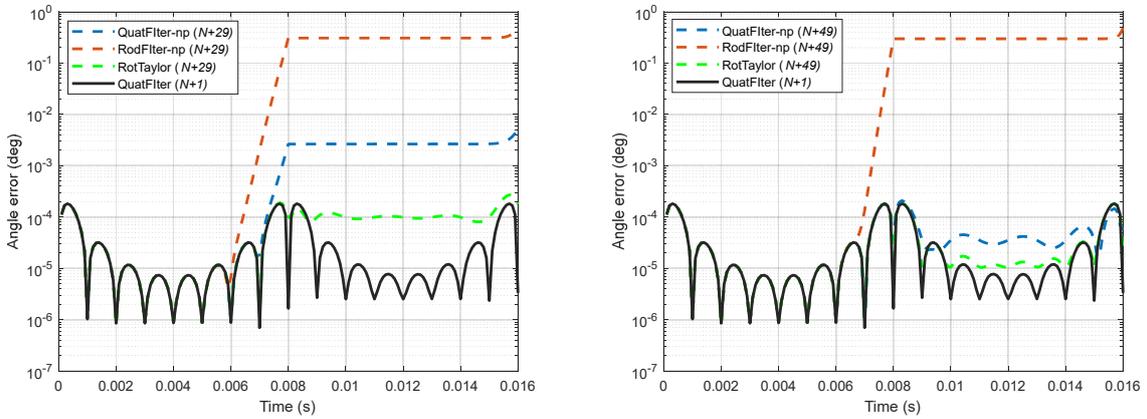

Figure 10. Attitude errors in two update intervals at coning frequency of 100 Hz for $N = 8$, with truncation order $m_T = N + 29$ (left) and $m_T = N + 49$ (right).

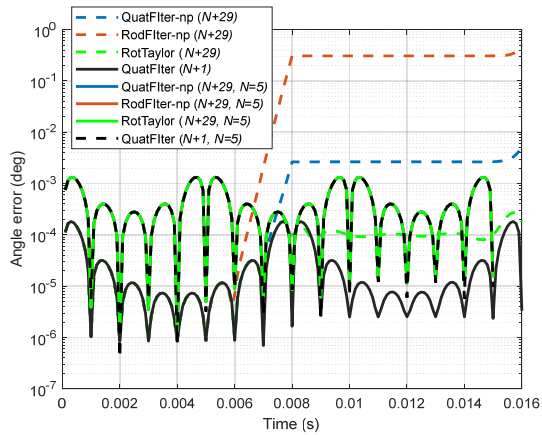

Figure 11. Attitude errors in three update intervals for $N = 5$ samples at coning frequency of 100 Hz.

RodTaylor, RodFIter and RotFIter are omitted in Fig. 9, because they are found to be nearly identical to those of QuatFIter-np, RodFIter-np and QuatFIter, respectively. We see that all algorithms have comparable accuracy while the coning frequency stays below 30 Hz, linearly increasing with respect to the relative frequency. Specifically, they all reach the machine precision for coning frequency less than 10 Hz. If the coning frequency goes up further, however, three algorithms of QuatFIter-np, RodFIter-np and RotTaylor, with the truncation order $m_T = N + 9$, begin to deteriorate and become even worse than the mainstream 2/3-sample algorithms for over 70 Hz. When the truncation order is set to $m_T = N + 29$, their accuracies do not decline until the coning frequency is larger than 70-80 Hz. In contrast, with a constant truncation order $m_T = N + 1$, the accuracy of QuatFIter (RodFIter and RotFIter) is uniformly the best over the whole frequency range, approaching those of the mainstream algorithms at the right end. It is observed in our simulations that higher truncation order would not bring further accuracy improvement to the QuatFIter (RodFIter and RotFIter). This significant advantage is largely owed to the excellent functional representation capability and numerical stability of the Chebyshev polynomial.

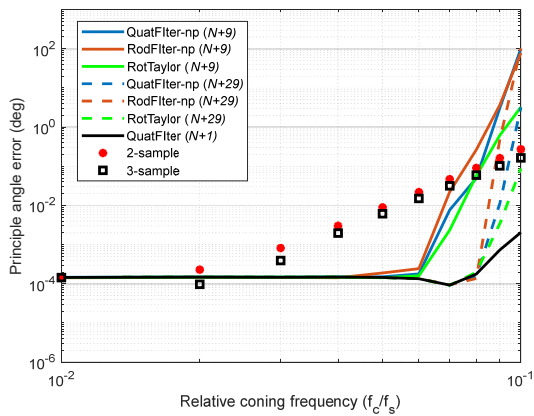

Figure 12. Attitude errors as function of relative frequency ($N = 8$) for navigation-grade gyroscope measurements.

The error rising of QuatFIter-np, RodFIter-np and RotTaylor appears odd, so the case of $f_c = 100$ Hz is particularly examined. Figure 10 presents their attitude errors for two update intervals using the polynomial truncation orders $m_T = N + 29$ and $m_T = N + 49$. With the truncation order increased, QuatFIter-np and RotTaylor improve in accuracy as expected but RodFIter-np changes little. Contrary to the prediction that their accuracies could be unlimitedly improved by further increasing the truncation order, we have observed that the three algorithms all encounter numerical failures when the truncation order is larger than 150. Of special interest is the 'u-shape' profile in the QuatFIter that has oscillated peaks close to both ends of iteration intervals. It is the famous Runge's phenomenon [54] that is ubiquitous in high-order polynomial interpolation for evenly-spaced samples, which can also be apparently identified in Fig. 2 for the fitted angular velocity with $N = 8$. An interesting thing is observed in the case of $N = 5$ samples, as shown in Fig. 11 (with the results of Fig. 10 as the background) in which the results of all algorithms overlap. In specific, QuatFIter-np and RodFIter-np demonstrate better accuracy than they do in Fig. 10 ($N = 8$). This unusual observation is believed to be incurred by the Runge's phenomenon, so is the numerical failure encountered above. Supposedly, the technique of depressing the Runge's phenomenon (e.g. using multiple integrals of gyroscope measurements [37, 51]) could be used to improve all derived algorithms including the already well-performing QuatFIter.

Finally, a practical situation with noisy gyroscope measurement is investigated, as the high-order/sample algorithms tend to be much more sensitive to narrow-band noises that might lead to pseudo-coning [37]. Noise errors with an angle random walk of $0.001 \ \deg/\sqrt{h}$, comparable to a navigation-grade inertial navigation system, are considered. A common set of random gyroscope noises are generated and fed to all algorithms for uniform comparison. Figure 12 plots the attitude errors in ten seconds as the function of relative frequency for $N = 8$. For coning frequencies below 20 Hz, all derived algorithms have similar accuracy with the mainstream 2/3-sample algorithms, as the noise dominates the attitude accuracy (cf. Fig. 9). The derived algorithms are comparable to each other when the coning frequency is below 40 Hz. The performance ranking keeps the same with the noise-free case in Fig. 9, except that the frequency point where QuatFIter-np, RodFIter-np and RotTaylor begin to get worse than QuatFIter increases from 30 Hz to 40 Hz.

## VIII. CONCLUSIONS AND DISCUSSIONS

This paper poses the strapdown attitude computation as seeking the general solutions to the kinematic equations of attitude parameters. Two basic approaches are briefly reviewed, namely, the Taylor series expansion approach and the functional iterative function approach (or alternatively known as the Picard iteration in the community of differential equations). Then, three groups of attitude algorithms have been derived by the two basic approaches, based on major attitude parameters including quaternion, Rodrigues vector and rotation vector. The first group, based on the Taylor series expansion (QuatTaylor, RodTaylor and RotTaylor), follows and considerably extends the framework of the mainstream algorithms, by making use of recursive calculation of high-order derivatives. The other two groups (QuatFIter-np, RodFIter-np and RotFIter-np; QuatFIter, RodFIter and RotFIter) both employ the functional iterative function approach yet with two different kinds of polynomial approximation, namely, the normal polynomial (presented first in this paper) and Chebyshev polynomial (recently published). Numerical tests under classical coning motions are carried out to compare the algorithms, refining the conclusions drawn in previous papers on the functional iterative integration approach. In the relative frequency range when the coning to sampling frequency ratio is below 0.05-0.1 (depending on the chosen polynomial truncation order), the three algorithm groups have the same order of accuracy if the same number of samples are used to fit the angular velocity over the iteration interval; in the range of higher relative frequency, the third group (Quat/Rod/RotFIter) performs better in both accuracy and robustness to the Runge phenomenon than the other two groups do, thanks to the unique properties of Chebyshev polynomial. Notably, the third group allows a lower truncation order, while the other two groups require significantly higher truncation order and might even encounter numerical failure.

All presented algorithms are in nature iterative and thus much computation-expensive relative to the mainstream algorithms. If necessary, however, they can be implemented in inertial navigation systems by optimized software or customized hardware for potential accuracy benefit. Additionally, all proposed algorithms are founded on polynomial approximation of angular velocity/specific force, so any potential improvement in polynomial approximation would further advance them, e.g., using the multiple integrals of gyroscopes/accelerometers measurements or resorting to the technique of depressing the Runge's phenomenon.

### Appendix: Mainstream 2/3-sample Algorithms

For the special case of $N = 2$ and $n = 1$, $\boldsymbol{\omega}(t) = \mathbf{d}_0 + \mathbf{d}_1 t$ according to (6), and the derivatives of the approximate rotation vector in (29) are explicitly

$$\begin{aligned}
\boldsymbol{\sigma}^{(1)}(0) &\approx \boldsymbol{\omega}(0) + \frac{1}{2}\boldsymbol{\omega}^{(-1)}(0) \times \boldsymbol{\omega}(0) = \mathbf{d}_0 \\
\boldsymbol{\sigma}^{(2)}(0) &\approx \boldsymbol{\omega}^{(1)}(0) + \frac{1}{2}\boldsymbol{\omega}^{(-1)}(0) \times \boldsymbol{\omega}^{(1)}(0) = \mathbf{d}_1 \\
\boldsymbol{\sigma}^{(3)}(0) &\approx \frac{1}{2}\boldsymbol{\omega}(0) \times \boldsymbol{\omega}^{(1)}(0) = \frac{1}{2}\mathbf{d}_0 \times \mathbf{d}_1 \\
\boldsymbol{\sigma}^{(j)}(0) &= 0, \ j \geq 4
\end{aligned} \tag{47}$$

Note that the notation $\boldsymbol{\omega}^{(-1)}(0) = \int_0^t \boldsymbol{\omega}\, dt \big|_{t=0} = 0$. Substituting into (3) yields the Taylor series approximation of the rotation vector

$$\boldsymbol{\sigma}(t) = \mathbf{d}_0 t + \frac{\mathbf{d}_1}{2} t^2 + \frac{\mathbf{d}_0 \times \mathbf{d}_1}{12} t^3 \tag{48}$$

Similarly, for the case of $N = 3$ and $n = 2$, $\boldsymbol{\omega}(t) = \mathbf{d}_0 + \mathbf{d}_1 t + \mathbf{d}_2 t^2$, and the derivatives of the approximate rotation vector in (29) are explicitly

$$\begin{aligned}
\boldsymbol{\sigma}^{(1)}(0) &\approx \mathbf{d}_0, \ \boldsymbol{\sigma}^{(2)}(0) \approx \mathbf{d}_1 \\
\boldsymbol{\sigma}^{(3)}(0) &\approx \boldsymbol{\omega}^{(2)}(0) + \frac{1}{2}\boldsymbol{\omega}(0) \times \boldsymbol{\omega}^{(1)}(0) = 2\mathbf{d}_2 + \frac{1}{2}\mathbf{d}_0 \times \mathbf{d}_1 \\
\boldsymbol{\sigma}^{(4)}(0) &\approx \boldsymbol{\omega}(0) \times \boldsymbol{\omega}^{(2)}(0) = 2\mathbf{d}_0 \times \mathbf{d}_2 \\
\boldsymbol{\sigma}^{(5)}(0) &\approx \boldsymbol{\omega}^{(1)}(0) \times \boldsymbol{\omega}^{(2)}(0) = 2\mathbf{d}_1 \times \mathbf{d}_2 \\
\boldsymbol{\sigma}^{(j)}(0) &= 0, \ j \geq 6
\end{aligned} \tag{49}$$

Substituting into (3) yields the Taylor series approximation of the rotation vector

$$\boldsymbol{\sigma}(t) = \mathbf{d}_0 t + \frac{\mathbf{d}_1}{2} t^2 + \frac{4\mathbf{d}_2 + \mathbf{d}_0 \times \mathbf{d}_1}{12} t^3 + \frac{\mathbf{d}_0 \times \mathbf{d}_2}{12} t^4 + \frac{\mathbf{d}_1 \times \mathbf{d}_2}{60} t^5 \tag{50}$$

Alternatively, the results in (48) and (50) can be directly obtained by integrating the second equation of (27).

When the angular increment measurements are considered, it can be obtained from (8) for $N = 2$ and 3, respectively,

$$\begin{bmatrix} \mathbf{d}_0 \\ \mathbf{d}_1 \end{bmatrix} = \begin{bmatrix} \left(3\Delta\boldsymbol{\theta}_1 - \Delta\boldsymbol{\theta}_2\right)\big/2T \\ \left(-\Delta\boldsymbol{\theta}_1 + \Delta\boldsymbol{\theta}_2\right)\big/T^2 \end{bmatrix} \quad \text{and} \quad \begin{bmatrix} \mathbf{d}_0 \\ \mathbf{d}_1 \\ \mathbf{d}_2 \end{bmatrix} = \begin{bmatrix} \left(11\Delta\boldsymbol{\theta}_1 - 7\Delta\boldsymbol{\theta}_2 + 2\Delta\boldsymbol{\theta}_3\right)\big/6T \\ \left(-2\Delta\boldsymbol{\theta}_1 + 3\Delta\boldsymbol{\theta}_2 - \Delta\boldsymbol{\theta}_3\right)\big/T^2 \\ \left(\Delta\boldsymbol{\theta}_1 - 2\Delta\boldsymbol{\theta}_2 + \Delta\boldsymbol{\theta}_3\right)\big/2T^3 \end{bmatrix} \tag{51}$$

Substituting (51) into (48) and (50) and letting $t = NT$ produce the well-known mainstream 2-sample and 3-sample algorithms [10], respectively,

$$\begin{aligned} N = 2: \quad & \boldsymbol{\sigma}(2T) = \Delta\boldsymbol{\theta}_1 + \Delta\boldsymbol{\theta}_2 + \frac{2}{3}\Delta\boldsymbol{\theta}_1 \times \Delta\boldsymbol{\theta}_2 \\ N = 3: \quad & \boldsymbol{\sigma}(3T) = \Delta\boldsymbol{\theta}_1 + \Delta\boldsymbol{\theta}_2 + \Delta\boldsymbol{\theta}_3 + 0.4125\Delta\boldsymbol{\theta}_1 \times \Delta\boldsymbol{\theta}_3 + 0.7125\Delta\boldsymbol{\theta}_2 \times \left(\Delta\boldsymbol{\theta}_3 - \Delta\boldsymbol{\theta}_1\right) \end{aligned} \tag{52}$$

ACKNOWLEDGEMENTS

Thanks to Dr. Qi Cai for Matlab code optimization.

Table I. Attitude Algorithms Derived by Taylor Series Expansion and Functional Iterative Integration

| | **Taylor Series Expansion** | **Functional Iterative Integration** | |
| --- | --- | --- | --- |
| | | *Normal Polynomial* | *Chebyshev Polynomial* (cf. [16, 33, 34]) |
| **Quaternion** | QuatTaylor<br>$\mathbf{q}(t)=\sum_{j=0}^{\infty}\mathbf{q}^{(j)}(0)\frac{t^j}{j!}$<br>$\mathbf{q}^{(j)}(0)=\frac{1}{2}\sum_{i=0}^{j-1}\binom{i-1}{j}\mathbf{q}^{(j-1-i)}(0)\circ\boldsymbol{\omega}^{(i)}(0)$ | QuatFIter-np<br>$\mathbf{q}_j(t)=\mathbf{q}(0)+\frac{1}{2}\sum_{k=0}^{m_{j-1}}\sum_{i=0}^{n}\frac{\mathbf{b}_{j-1,k}\circ\mathbf{d}_i}{k+i+1}t^{k+i+1}$ | QuatFIter<br>$\mathbf{q}_j(\tau)=\mathbf{q}(0)+\frac{t_N}{8}\sum_{k=0}^{m_{j-1}}\sum_{i=0}^{n}\mathbf{b}_{j-1,k}\circ\mathbf{c}_i\left(G_{k+i[-1\,\tau]}+G_{|k-i|[-1\,\tau]}\right)$ |
| **Rodrigues Vector** | RodTaylor<br>$\mathbf{g}(t)=\sum_{j=0}^{\infty}\mathbf{g}^{(j)}(0)\frac{t^j}{j!}$<br>$\mathbf{g}^{(j)}(0)=\boldsymbol{\omega}^{(j-1)}(0)+\frac{1}{2}\sum_{i=0}^{j-1}\binom{i-1}{j}\mathbf{g}^{(j-1-i)}(0)\times\boldsymbol{\omega}^{(i)}(0)$<br>$+\frac{1}{4}\sum_{i=0}^{j-1}\sum_{k=0}^{i}\binom{i-1}{j}\binom{i}{k}\mathbf{g}^{(j-1-i)}(0)\left(\mathbf{g}^{(i-k)}(0)\right)^T\boldsymbol{\omega}^{(k)}(0)$ | RodFIter-np<br>$\mathbf{g}_j(t)=\sum_{i=0}^{n}\frac{\mathbf{d}_i}{i+1}t^{i+1}+\frac{1}{2}\sum_{k=0}^{m_{j-1}}\sum_{i=0}^{n}\frac{\mathbf{b}_{j-1,k}\times\mathbf{d}_i}{k+i+1}t^{k+i+1}$<br>$+\frac{1}{4}\sum_{s=0}^{m_{j-1}}\sum_{k=0}^{m_{j-1}}\sum_{i=0}^{n}\frac{\mathbf{b}_{j-1,s}^T\mathbf{d}_i}{s+k+i+1}t^{s+k+i+1}$ | RodFIter<br>$\mathbf{g}_j(\tau)=\frac{t_N}{2}\left(\begin{array}{l}\sum_{i=0}^{n}\mathbf{c}_iG_{i[-1\,\tau]}+\frac{1}{4}\sum_{k=0}^{m_{j-1}}\sum_{i=0}^{n}\mathbf{b}_{j-1,k}\times\mathbf{c}_i\left(G_{k+i[-1\,\tau]}+G_{|k-i|[-1\,\tau]}\right)\\+\frac{1}{16}\sum_{s=0}^{m_{j-1}}\sum_{k=0}^{m_{j-1}}\sum_{i=0}^{n}\mathbf{b}_{j-1,s}\mathbf{b}_{j-1,k}^T\mathbf{c}_i\left(\begin{array}{l}G_{s+k+i[-1\,\tau]}+G_{|s+k-i|[-1\,\tau]}\\+G_{|k-s|+i[-1\,\tau]}+G_{|k-s|-i|[-1\,\tau]}\end{array}\right)\end{array}\right)$ |
| **Rotation Vector** | RotTaylor<br>$\boldsymbol{\sigma}(t)=\sum_{j=0}^{\infty}\boldsymbol{\sigma}^{(j)}(0)\frac{t^j}{j!}$<br>$\boldsymbol{\sigma}^{(j)}(0)=\boldsymbol{\omega}^{(j-1)}(0)+\frac{1}{2}\sum_{i=0}^{j-1}\binom{i-1}{j}\boldsymbol{\sigma}^{(j-1-i)}(0)\times\boldsymbol{\omega}^{(i)}(0)$<br>$+\sum_{i=0}^{j-1}\sum_{k=0}^{i}\sum_{s=0}^{k}\binom{i-1}{j}\binom{i}{k}\binom{k}{s}A^{(j-1-i)}(0)\boldsymbol{\sigma}^{(i-k)}(0)\times\left(\boldsymbol{\sigma}^{(k-s)}(0)\times\boldsymbol{\omega}^{(s)}(0)\right)$<br><br>RotTaylor-T2/T2s<br>$\boldsymbol{\sigma}^{(j)}(0)=\boldsymbol{\omega}^{(j-1)}(0)+\frac{1}{2}\sum_{i=0}^{j-1}\binom{i-1}{j}\boldsymbol{\sigma}^{(j-1-i)}(0)\times\boldsymbol{\omega}^{(i)}(0)$<br>$\boldsymbol{\sigma}^{(j)}(0)=\boldsymbol{\omega}^{(j-1)}(0)+\frac{1}{2}\sum_{i=0}^{j-1}\binom{i-1}{j}\boldsymbol{\omega}^{(j-2-i)}(0)\times\boldsymbol{\omega}^{(i)}(0)$ | RotFIter-np-T3/T2<br>$\boldsymbol{\sigma}_j(t)=\sum_{i=0}^{n}\frac{\mathbf{d}_i}{i+1}t^{i+1}+\frac{1}{2}\sum_{k=0}^{m_{j-1}}\sum_{i=0}^{n}\frac{\mathbf{b}_{j-1,k}\times\mathbf{d}_i}{k+i+1}t^{k+i+1}$<br>$+\frac{1}{12}\sum_{s=0}^{m_{j-1}}\sum_{k=0}^{m_{j-1}}\sum_{i=0}^{n}\frac{\mathbf{b}_{j-1,s}\times(\mathbf{b}_{j-1,k}\times\mathbf{d}_i)}{s+k+i+1}t^{s+k+i+1}$<br>$\boldsymbol{\sigma}_j(t)=\sum_{i=0}^{n}\frac{\mathbf{d}_i}{i+1}t^{i+1}+\frac{1}{2}\sum_{k=0}^{m_{j-1}}\sum_{i=0}^{n}\frac{\mathbf{b}_{j-1,k}\times\mathbf{d}_i}{k+i+1}t^{k+i+1}$ | RotFIter<br>$\boldsymbol{\sigma}_j(\tau)=\frac{t_N}{2}\left(\begin{array}{l}\sum_{i=0}^{n}\mathbf{c}_iG_{i[-1\,\tau]}+\sum_{k=0}^{p_{j-1}}\boldsymbol{\gamma}_{j-1,k}G_{k[-1\,\tau]}\\+\frac{1}{4}\sum_{k=0}^{m_{j-1}}\sum_{i=0}^{n}\mathbf{b}_{j-1,k}\times\mathbf{c}_i\left(G_{k+i[-1\,\tau]}+G_{|k-i|[-1\,\tau]}\right)\end{array}\right)$<br>$\boldsymbol{\gamma}_{j-1,k}\approx\frac{2-\delta_{0k}}{Q}\sum_{s=0}^{Q-1}\cos\left(\frac{k(s+1/2)\pi}{Q}\right)\boldsymbol{\eta}_{j-1}\left(\cos\left(\frac{(s+1/2)\pi}{Q}\right)\right)$<br>$\boldsymbol{\eta}_{j-1}\triangleq A\boldsymbol{\sigma}_{j-1}\times(\boldsymbol{\sigma}_{j-1}\times\boldsymbol{\omega})$<br><br>RotFIter-T3/T2<br>$\boldsymbol{\sigma}_j(\tau)=\frac{t_N}{2}\left(\begin{array}{l}\sum_{i=0}^{n}\mathbf{c}_iG_{i[-1\,\tau]}+\frac{1}{4}\sum_{k=0}^{m_{j-1}}\sum_{i=0}^{n}\mathbf{b}_{j-1,k}\times\mathbf{c}_i\left(G_{k+i[-1\,\tau]}+G_{|k-i|[-1\,\tau]}\right)\\+\frac{1}{48}\sum_{s=0}^{m_{j-1}}\sum_{k=0}^{m_{j-1}}\sum_{i=0}^{n}\mathbf{b}_{j-1,s}\times\mathbf{b}_{j-1,k}\times\mathbf{c}_i\left(\begin{array}{l}G_{s+k+i[-1\,\tau]}+G_{|s+k-i|[-1\,\tau]}\\+G_{|k-s|+i[-1\,\tau]}+G_{|k-s|-i|[-1\,\tau]}\end{array}\right)\end{array}\right)$<br>$\boldsymbol{\sigma}_j(\tau)=\frac{t_N}{2}\left(\sum_{i=0}^{n}\mathbf{c}_iG_{i[-1\,\tau]}+\frac{1}{4}\sum_{k=0}^{m_{j-1}}\sum_{i=0}^{n}\mathbf{b}_{j-1,k}\times\mathbf{c}_i\left(G_{k+i[-1\,\tau]}+G_{|k-i|[-1\,\tau]}\right)\right)$ |